\documentclass[11pt]{amsart}
\usepackage{amsmath}
\usepackage{amssymb}
\usepackage{amsthm} 
\usepackage{bbm}
\usepackage{graphicx}
\usepackage{algorithm}
\usepackage{algorithmic}
\usepackage{clrscode}
\usepackage{fullpage}
\usepackage{stmaryrd}
\usepackage{color}
\usepackage{url}
\usepackage{hyperref}
\usepackage{caption}
\usepackage{subcaption}

\newtheorem{prop}{Proposition}

\newtheorem{thm}{Theorem}

\newtheorem{Def}{Definition}

\def\x{{\mathbf x}}

\DeclareMathOperator*{\argmin}{arg\,min}

\begin{document}
\pagestyle{plain}
\newenvironment{frcseries}{\fontfamily{frc} \selectfont}{}
\newcommand{\textfrc}[1]{{\frcseries #1}}
\newcommand{\mathfrc}[1]{\text{\textfrc{#1}}}

% Title.
% ------
\title{Robust Sparse Phase Retrieval Made Easy}\thanks{M.A. Iwen:  Department of Mathematics and Department of ECE, Michigan State University ({\tt markiwen@math.msu.edu}).  M.A. Iwen was supported in part by NSF DMS-1416752 and NSA H98230-13-1-0275.\\  \indent A. Viswanathan:  Department of Mathematics, Michigan State University ({\tt aditya@math.msu.edu}). \\  \indent Y. Wang: Department of Mathematics, The Hong Kong University of Science and Technology ({\tt yangwang@ust.hk}). Y. Wang was partially supported by NSF DMS-1043032 and AFOSR FA9550-12-1-0455.}
%
% Single address.
% ---------------
\author{
{Mark Iwen} %
\and Aditya Viswanathan
\and Yang Wang
}

\begin{abstract}
In this short note we propose a simple two-stage sparse phase retrieval strategy that uses a near-optimal number of measurements, and is both computationally efficient and robust to measurement noise.  In addition, the proposed strategy is fairly general, allowing for a large number of new measurement constructions and recovery algorithms to be designed with minimal effort.
\end{abstract}

\maketitle
\thispagestyle{empty}

\section{Introduction}
\label{sec:Intro}

Herein we consider the phase retrieval problem of reconstructing a given vector ${\bf x} \in \mathbbm{C}^N$ from noisy magnitude measurements of the form 
\begin{equation}
b_i := \left| \langle {\bf p}_i, {\bf x}  \rangle \right|^2 + n_i,
\label{equ:PRmeasureIntro}
\end{equation}
where ${\bf p}_i \in \mathbbm{C}^N$ is a measurement vector, and $n_i \in \mathbbm{R}$ represents arbitrary measurement noise, for $i = 1, \dots, M$.  In particular, we focus on the setting where the dimension $N$ is either very large, or else the number of measurements allowed, $M$, is otherwise severely restricted.  In either case, our inability to gather the $M = \mathcal{O}(N)$ measurements required for the recovery of ${\bf x}$ in general \cite{heinosaari2013quantum} forces us to consider the possibility of approximating ${\bf x}$ using only $M \ll N$ magnitude measurements, if possible.  This is the situation motivating the {\em compressive phase retrieval problem} (see, e.g., \cite{ohlsson2012cprl,SCHNITNER!,li2013sparse,jaganathan2013sparse,wang2014phase,eldar2014phase,GESPAR2014Eldar,yapar2014fast}), in which one attempts to accurately approximate ${\bf x} \in \mathbbm{C}^N$ using only $M = o(N)$ magnitude measurements $\eqref{equ:PRmeasureIntro}$ under the assumption that ${\bf x}$ is either sparse, or compressible.

One question regarding the compressive phase retrieval problem is how many measurements are needed to allow for stable reconstruction of $\x$.   Clearly, compressive phase retrieval requires at least as many measurements as the corresponding classical compressive sensing problem since one is given less information.  Hence, stable compressive phase retrieval requires at least $\mathcal{O}(s\log(N/s))$ magnitude measurements\footnote{See, e.g., Chapter 10 of \cite{HolgerBook} concerning the minimal number of measurements required for stable compressive sensing.} -- but {\it can it be done with $M=\mathcal{O}(s\log(N/s))$ measurements}?  It is shown in \cite{eldar2014phase} that stable compressive phase retrieval is indeed achievable with $M=\mathcal{O}(s\log(N/s))$ measurements for {\it real} $\x$ if the entries of ${\bf p}_i$ are real independent and identically distributed (i.i.d.) Gaussians.  However, this question was unresolved in the complex case.  In this note we extend the result to the complex case.  Furthermore, we do so in a constructive way by providing a computational procedure which can stably reconstruct complex $\x$ using only $\mathcal{O}(s\log(N/s))$ magnitude measurements.

Unlike previous sparse phase retrieval approaches, we propose a generic two-stage solution technique consisting of (i) using the phase retrieval technique of one's choice to recover compressive sensing measurements of ${\bf x}$, $\mathcal{C}{\bf x}  \in \mathbbm{C}^m$, followed by (ii) utilizing the compressive sensing method of one's choice in order to approximate ${\bf x}$ from the recovered measurements $\mathcal{C}{\bf x}$.  As we shall see, the generic nature of the proposed sparse phase retrieval procedure not only allows for a relatively large number of measurement matrices and recovery algorithms to be used, but also allows robust recovery guarantees for the sparse phase retrieval problem to be proven in the complex setting essentially ``for free'' by combining existing robust recovery results from the compressive sensing literature with robust recovery results for the standard phase retrieval setting.  As a result, we are able to show that $\mathcal{O}(s \log(N/s))$ magnitude measurements suffice in order to recover a large class of compressible vectors with the same quality of error guarantee as commonly achieved in the compressive sensing literature.  Finally, numerical experiments demonstrate that the proposed approach is also both efficient and robust in practice.

\section{Background}
\label{sec:background}

In this section we briefly recall selected results from the existing literature on compressive sensing \cite{donohoCS, HolgerBook} and phase retrieval \cite{balan2006signal,balan2009painless,candes2013phaselift,candes2014solving,alexeev2014phase,FMNW2014}. Let $\| {\bf x} \|_0$ denote the number of nonzero entries in a given ${\bf x } \in \mathbbm{C}^N$, and $\| {\bf x} \|_p$ denote the standard $\ell_p$-norm of ${\bf x }$ for all $p \geq 1$, i.e., $\| {\bf x} \|_p := \left(\sum^{N}_{n=1} |x_n|^p \right)^{1/p}$ for all ${\bf x} \in \mathbbm{C}^N$.

\subsection{Compressive Sensing}

Compressive sensing methods deal with the construction of an $m \times N$ measurement matrix, $\mathcal{C}$, with $m$ minimized as much as possible subject to the constraint that an associated approximation algorithm, $\Delta_{\mathcal{C}}: \mathbbm{C}^m \rightarrow \mathbbm{C}^N$, can still accurately approximate any given vector ${\bf x} \in \mathbbm{C}^N$.  More precisely, compressive sensing methods allow one to minimize $m$, the number of rows in $\mathcal{C}$, as a function of $s$ and $N$ such that
\begin{equation}
\left\| ~ \Delta_{\mathcal{C}} \left( \mathcal{C} {\bf x} \right) - {\bf x} ~ \right\|_p \leq C_{p,q} \cdot s^{\frac{1}{p} - \frac{1}{q}} \left( \inf_{{\bf z} \in \mathbbm{C}^N, \| {\bf z} \|_0 \leq s} \left\| {\bf x} - {\bf z} \right\|_q \right)
\label{eqn:Aerror}
\end{equation}
holds for all ${\bf x} \in \mathbbm{C}^N$ in various fixed $\ell_p$,$\ell_q$ norms, $1 \leq q \leq p \leq 2$, for an absolute constant $C_{p,q} \in \mathbbm{R}$ (e.g., see \cite{BestkTerm, HolgerBook}).  Note that this implies that $ {\bf x}$ will be recovered exactly if it contains only $s$ nonzero entries.  Similarly, $ {\bf x}$ will be accurately approximated by $\Delta_{\mathcal{C}} \left( \mathcal{C}  {\bf x} \right)$ any time its $\ell_q$-norm is dominated by its largest $s$ entries.

There are a wide variety of measurement matrices $\mathcal{C} \in \mathbbm{C}^{m \times N}$ with $m = \mathcal{O}(s \log (N/s))$ that have associated approximation algorithms, $\Delta_{\mathcal{C}}$, which are computationally efficient, numerically robust, and able to achieve error guarantees of the form \eqref{eqn:Aerror} for all ${\bf x} \in \mathbbm{C}^N$.  For example, this is true of ``most'' random matrices $\mathcal{C} \in \mathbbm{C}^{m \times N}$ with i.i.d. subgaussian random entries \cite{baraniuk2008simple,HolgerBook}.  Similarly, one may construct such a $\mathcal{C} \in \mathbbm{C}^{m \times N}$ with high probability by selecting a set of $m = \mathcal{O}(s \log^4 N)$ rows uniformly at random from an $N \times N$ discrete Fourier transform matrix (or, more generally, from any ``sufficiently flat'' $N \times N$ unitary matrix) \cite{HolgerBook}.  In either case, one may then use a large number of approximation algorithms, $\Delta_{\mathcal{C}}$, that will achieve error guarantees along the lines of \eqref{eqn:Aerror}, including convex optimization techniques \cite{CS1, CS4, NearOpt}, iterative hard thresholding \cite{HardThreshforCS}, (regularized) orthogonal matching pursuit \cite{CS2,OMPvsBP,ROMP,ROMPstable}, and the CoSaMP algorithm \cite{COSAMP}, to name just a few.

More generally, any matrix with the \textit{robust null space property} \cite{BestkTerm} will have an associated approximation algorithm that is both computationally efficient and numerically robust.  Let $S\subseteq\{1,2,...,N\}$, and ${\bf x} \in \mathbbm{C}^N$.  Then, ${\bf x}_{S}$ will denote ${\bf x}$ with all entries not in $S$ set to zero.  That is,
$$
(x_{S})_{j} :=
\begin{cases}
0, & \text{if } j\notin S, \\
x_j, & \text{if } j\in S.\\
\end{cases}
$$
The robust null space property can now be defined as follows.

\begin{Def}
Let $s, m, N \in \mathbbm{N}$ be such that $s < m < N$.  We will say that the matrix $\mathcal{C} \in \mathbbm{C}^{m \times N}$ satisfies the $\ell_2$-robust null space property of order $s$ with constants $0 < \rho < 1$ and $\tau > 0$ if 
$$\| {\bf x}_S \|_1 \leq \rho \| {\bf x}_{S^c} \|_1 + \tau \| \mathcal{C} {\bf x} \|_2$$
holds for all ${\bf x} \in \mathbbm{C}^{N}$ and $S \subset \{ 1,2,...,N \}$ with cardinality $|S| \leq s$, where $S^c$ denotes the complement of $S$.  
\label{def:rnsp}
\end{Def}
In particular, the following robust compressive sensing result for matrices with the null space property is a restatement of Theorem 4.22 from \cite{HolgerBook}.

\begin{thm}
Suppose that the matrix $\mathcal{C} \in \mathbbm{C}^{m \times N}$ satisfies the $\ell_2$-robust null space property of order $s$ with constants $0 < \rho < 1$ and $\tau > 0$.  Then, for any ${\bf x} \in \mathbbm{C}^N$, the vector 
\begin{equation}
{\bf \tilde{x}} := \argmin_{{\bf z} \in \mathbbm{C}^N}~ \| {\bf z} \|_1 ~\textrm{~subject to~}~ \| \mathcal{C}{\bf z} - {\bf y} \|_2 \leq \eta,
\label{equ:xopt}
\end{equation}
where ${\bf y} := \mathcal{C}{\bf x} + {\bf e}$ for some ${\bf e} \in \mathbbm{C}^m$ with $\| {\bf e} \|_2 \leq \eta$, will satisfy
\begin{equation}
\| {\bf x} - {\bf \tilde{x}} \|_2 \leq \frac{C}{\sqrt{s}} \cdot \left( \inf_{{\bf z} \in \mathbbm{C}^N, \| {\bf z} \|_0 \leq s} \left\| {\bf x} - {\bf z} \right\|_1 \right) + D \eta
\end{equation}
for some constants $C,D \in \mathbbm{R}^+$ that only depend on $\rho$ and $\tau$. 
\label{thm:CS}
\end{thm}

Many matrices exist with the $\ell_2$-robust null space property including, e.g., ``most'' randomly constructed subgaussian and subsampled discrete Fourier transform matrices (as per above).  Thus, in some sense it is not difficult to find a matrix $\mathcal{C} \in \mathbbm{C}^{m \times N}$ to which Thoerem~\ref{thm:CS} will apply.  Furthermore, ${\bf \tilde{x}}$ from \eqref{equ:xopt} can be computed efficiently via convex optimization techniques.  See \cite{HolgerBook} for details.

\subsection{Phase Retrieval}

Noisy phase retrieval problems involve the reconstruction of a given vector ${\bf x} \in \mathbbm{C}^N$, up to a global phase factor, from magnitude measurements of the form
\begin{equation}
b_i := \left| \langle {\bf p}_i, {\bf x}  \rangle \right|^2 + n_i,
\label{equ:PRmeasure}
\end{equation}
where ${\bf p}_i \in \mathbbm{C}^N$ and $n_i \in \mathbbm{R}$ for $i = 1, \dots, M$.  Vectorizing \eqref{equ:PRmeasure} yields 
\begin{equation}
{\bf b} := \left| \mathcal{P} {\bf x} \right|^2 + {\bf n}, 
\label{equ:PRmeas}
\end{equation}
where ${\bf b},{\bf n} \in \mathbbm{R}^M$, $\mathcal{P} \in \mathbbm{C}^{M \times N}$, and $| \cdot |^2: \mathbbm{C}^M \rightarrow \mathbbm{R}^M$ computes the  component-wise squared magnitude of each vector entry.  Thus, the primary objective of phase retrieval is to construct a recovery algorithm, $\Phi_{\mathcal{P}}:  \mathbbm{R}^M \rightarrow \mathbbm{C}^N$, that satisfies a relative error guarantee such as, e.g.,
\begin{equation}
\min_{\theta \in [0, 2 \pi]} \left( \frac{\left\| ~ \Phi_{\mathcal{P}} \left( {\bf b} \right) - \mathbbm{e}^{\mathbbm{i} \theta} {\bf x} ~ \right\|_2}{ \| {\bf x} \|_2} \right)^q \leq C_{\mathcal{P}} \cdot \frac{\| {\bf n} \|_2 }{ \sqrt{M} \| {\bf x} \|^2_2}
\label{eqn:PRerror}
\end{equation}
for a particular measurement matrix $\mathcal{P} \in \mathbbm{C}^{M \times N}$, $q \in [1, 2]$, and approximation factor $C_{\mathcal{P}} \in \mathbbm{R}^+$ (which may depend on $\mathcal{P}$).  

Several recovery algorithms achieve error guarantees along the lines of \eqref{eqn:PRerror} while using at most $M = \mathcal{O}(N \log N)$ measurements, including both {\em PhaseLift} \cite{candes2013phaselift,candes2014solving} as well as a more recent graph-theoretic and frame-based approach \cite{alexeev2014phase}.  In particular, the following robust phase retrieval result is a variant of Theorem 1.3 from \cite{candes2014solving}.\footnote{Equation (1.8) in Theorem 1.3 is technically incorrect as stated in \cite{candes2014solving}.  See \cite{PhaseLiftCorrection} for a corrected and simplified proof of Theorem~\ref{thm:PR} as stated herein.}
\begin{thm}
Let $\mathcal{P} \in \mathbbm{C}^{M \times N}$ have its $M$ rows be independently drawn either uniformly at random from the sphere of radius $\sqrt{N}$ in $\mathbbm{C}^N$, or else as complex normal random vectors from $\mathcal{N}(0,\mathcal{I}_N / 2) + \mathbbm{i} \mathcal{N}(0,\mathcal{I}_N / 2)$.  Then, $\exists$ universal constants $\tilde{B}, \tilde{C}, \tilde{D} \in \mathbbm{R}^+$ such that the PhaseLift procedure $\Phi_{\mathcal{P}}:  \mathbbm{R}^M \rightarrow \mathbbm{C}^N$ satisfies 
\begin{equation}
\min_{\theta \in [0, 2 \pi]}\left\| ~ \Phi_{\mathcal{P}} \left( {\bf b} \right) - \mathbbm{e}^{\mathbbm{i} \theta} {\bf x} ~ \right\|_2 \leq \tilde{C} \cdot \frac{\| {\bf n} \|_1 }{M \| {\bf x} \|_2}
\end{equation}
for all ${\bf x} \in \mathbbm{C}^N$ with probability $1 - \mathcal{O}( \mathbbm{e}^{-\tilde{B} M} )$, provided that $M \geq \tilde{D} N$.  Here ${\bf b}, {\bf n} \in \mathbbm{R}^M$ are as in \eqref{equ:PRmeas}.
\label{thm:PR}
\end{thm}

Finally, it is important to note that the {\em PhaseLift} procedure from Theorem~\ref{thm:PR} can be computed via semidefinite programming techniques.  Thus, it is computationally tractable for modest dimensions, $N$.  See \cite{candes2013phaselift,candes2014solving} for details.

\section{A Simple Two-Stage Technique for Sparse Phase Retrieval}
\label{sec:Method}

In this section we consider using noisy magnitude measurements of the form
\begin{equation}
{\bf b} := \left| \mathcal{P} \mathcal{C} {\bf x} \right|^2 + {\bf n}, 
\label{equ:CPRmeas}
\end{equation}
where $\mathcal{P} \in \mathbbm{C}^{\tilde{m} \times m}$ is any phase retrieval matrix with an associated recovery algorithm $\Phi_{\mathcal{P}}:  \mathbbm{R}^{\tilde{m}} \rightarrow \mathbbm{C}^m$ that has an error guarantee along the lines of \eqref{eqn:PRerror}, and $\mathcal{C} \in \mathbbm{C}^{m \times N}$ is any compressive sensing matrix with an associated approximation algorithm $\Delta_{\mathcal{C}}: \mathbbm{C}^m \rightarrow \mathbbm{C}^N$ that has an error guarantee like \eqref{eqn:Aerror}.  In this situation the composition of the two recovery algorithms, $\Delta_{\mathcal{C}} \circ \Phi_{\mathcal{P}}:  \mathbbm{R}^{\tilde{m}} \rightarrow \mathbbm{C}^N$, should accurately approximate ${\bf x} \in \mathbbm{C}^N$, up to a global phase factor, from ${\bf b}$ whenever ${\bf x}$ is sufficiently sparse or compressible.  This leads us to the following intuitive observation.
\begin{prop}
Let ${\mathcal A} = \mathcal{P} \mathcal{C}$ where $\mathcal{C} \in \mathbbm{C}^{m\times N}$ has the robust null space property and $\mathcal{P} \in \mathbbm{C}^{\tilde{m}\times m}$ is a stable phase retrieval matrix. Then, $\mathcal{A}$ has the stable compressive phase retrieval property.
\label{Prop}
\end{prop}
More specifically, the following compressive phase retrieval result follows easily from Theorems~\ref{thm:CS} and~\ref{thm:PR}.

\begin{thm}
Let $\mathcal{P} \in \mathbbm{C}^{\tilde{m} \times m}$ have its $\tilde{m}$ rows be independently drawn either uniformly at random from the sphere of radius $\sqrt{m}$ in $\mathbbm{C}^m$, or else as complex normal random vectors from $\mathcal{N}(0,\mathcal{I}_m / 2) + \mathbbm{i} \mathcal{N}(0,\mathcal{I}_m / 2)$.  Furthermore, suppose that $\mathcal{C} \in \mathbbm{C}^{m \times N}$ satisfies the $\ell_2$-robust null space property of order $s$ with constants $0 < \rho < 1$ and $\tau > 0$.  Then, there exists a phase retrieval procedure, $\Phi_{\mathcal{P}}:  \mathbbm{R}^{\tilde{m}} \rightarrow \mathbbm{C}^m$, and a compressive sensing recovery algorithm, $\Delta_{\mathcal{C}}: \mathbbm{C}^m \rightarrow \mathbbm{C}^N$, such that
\begin{equation}
\min_{\theta \in [0, 2 \pi]} \left \| \mathbbm{e}^{\mathbbm{i} \theta} {\bf x} - \Delta_{\mathcal{C}} \left( \Phi_{\mathcal{P}} \left( {\bf b} \right) \right) \right \|_2 \leq \frac{C}{\sqrt{s}} \cdot \left( \inf_{{\bf z} \in \mathbbm{C}^N, \| {\bf z} \|_0 \leq s} \left\| {\bf x} - {\bf z} \right\|_1 \right) + D \cdot \frac{\| {\bf n} \|_1 }{\tilde{m} \| \mathcal{C}{\bf x} \|_2}
\label{ThmCPR:error}
\end{equation}
holds for all ${\bf x} \in \mathbbm{C}^N$ with probability $1 - \mathcal{O}( \mathbbm{e}^{-B \tilde{m}} )$, provided that $\tilde{m} \geq E \cdot m$.  Here ${\bf b}, {\bf n} \in \mathbbm{R}^{\tilde{m}}$ are as in \eqref{equ:CPRmeas}, and $B, E \in \mathbbm{R}^+$ are universal constants, while $C,D \in \mathbbm{R}^+$ are constants that only depend on $\rho$ and $\tau$. 
\label{thm:CPR}
\end{thm}

Considering the number of magnitude measurements required by Theorem~\ref{thm:CPR}, we note that $\tilde{m} = \mathcal{O}(s \log (N/s))$ such measurements will suffice to achieve \eqref{ThmCPR:error} for all ${\bf x} \in \mathbbm{C}^N$ with high probability whenever $\mathcal{C} \in \mathbbm{C}^{m \times N}$ is, e.g., a random matrix with i.i.d. subgaussian random entries.  In this situation $\mathcal{C}$ will also likely have both (i) the $\ell_2$-robust null space property of order $s$ with constants $0 < \rho < 1$ and $\tau > 0$, and (ii) a small \textit{restricted isometry constant of order $2s$}, $\delta_{2s} < 1$ (see, e.g., \S 6.2 and \S9.1 of \cite{HolgerBook} for details).  As a consequence, $\mathcal{C}$ will also satisfy
\begin{equation}
\frac{1}{\tau} \cdot \max_{S \subset \{ 1, \dots, N \}, |S| = s} \left( \| {\bf x}_S \|_1 -  \rho \| {\bf x}_{S^c} \|_1 \right) ~\leq~ \| \mathcal{C}{\bf x} \|_2 ~\leq~ \sqrt{1 - \delta_{2s}} \left( \| {\bf x} \|_2 + \frac{\| {\bf x} \|_1}{\sqrt{2s}} \right) 
\label{equ:Cxbounds}
\end{equation}
for all ${\bf x} \in \mathbbm{C}^N$ with high probability (w.h.p.).\footnote{The lower bound is a simple consequence of Definition~\ref{def:rnsp}.  For the upper bound see, e.g., Exercise 6.6 in \cite{HolgerBook}.}  Considering the Theorem~\ref{thm:CPR} error guarantee \eqref{ThmCPR:error} in light of \eqref{equ:Cxbounds}, we can now see that Theorem~\ref{thm:CPR} implies that all sufficiently compressible vectors with, e.g., 
\begin{equation}
\frac{1}{\sqrt{\tilde{m}}} ~\leq~ \frac{1}{\tau} \cdot \max_{S \subset \{ 1, \dots, N \}, |S| = s} \left( \| {\bf x}_S \|_1 -  \rho \| {\bf x}_{S^c} \|_1 \right)
\end{equation}
will also satisfy 
\begin{equation}
\min_{\theta \in [0, 2 \pi]} \left \| \mathbbm{e}^{\mathbbm{i} \theta} {\bf x} - \Delta_{\mathcal{C}} \left( \Phi_{\mathcal{P}} \left( {\bf b} \right) \right) \right \|_2 \leq \frac{C}{\sqrt{s}} \cdot \left( \inf_{{\bf z} \in \mathbbm{C}^N, \| {\bf z} \|_0 \leq s} \left\| {\bf x} - {\bf z} \right\|_1 \right) + D \| {\bf n} \|_2
\label{ThmCPR:errorSparse}
\end{equation}
w.h.p. whenever $\mathcal{C}$ is a random matrix with i.i.d. subgaussian entries.

Finally, it is interesting to note that the two-stage approach outlined in this section also confers some computational advantages.  Mainly, the phase retrieval recovery algorithm $\Phi_{\mathcal{P}}:  \mathbbm{R}^{\tilde{m}} \rightarrow \mathbbm{C}^m$ only needs to recover a vector of length $m = \mathcal{O}(s \log(N/s))$.  This allows phase retrieval approaches based on, e.g., semidefinite programming to efficiently approximate significantly larger vectors ${\bf x} \in \mathbbm{C}^N$ than otherwise possible when $N \gg s$.  

\section{Empirical Evaluation}
\label{sec:Eval}
We now present representative results demonstrating the numerical robustness 
and efficiency of the proposed two-step strategy. For the results in this section, 
we use {\em PhaseLift} \cite{candes2013phaselift,candes2014solving} and Basis 
Pursuit \cite{CS4} to solve the phase retrieval and compressive sensing problems in 
steps (i) and (ii), respectively. Moreover, we use complex Gaussian phase retrieval 
matrices $\mathcal P$ and real Gaussian compressive sensing matrices 
$\mathcal C$. Matlab code used to generate the numerical results -- implemented 
using the optimization software packages TFOCS \cite{tfocs,tfocspaper} and CVX 
\cite{cvx,cvxpaper} -- is freely available at \cite{bitbucket_sparsePR}.

In each of the following results, we recover sparse, unit-norm complex vectors 
whose non-zero indices are independently and randomly chosen, and, whose 
non-zero entries are  i.i.d. standard complex Gaussians.
\begin{figure}[hbtp]
\centering
\includegraphics[clip=true, trim = 1.25in 0in 0in 0in, scale=0.45]{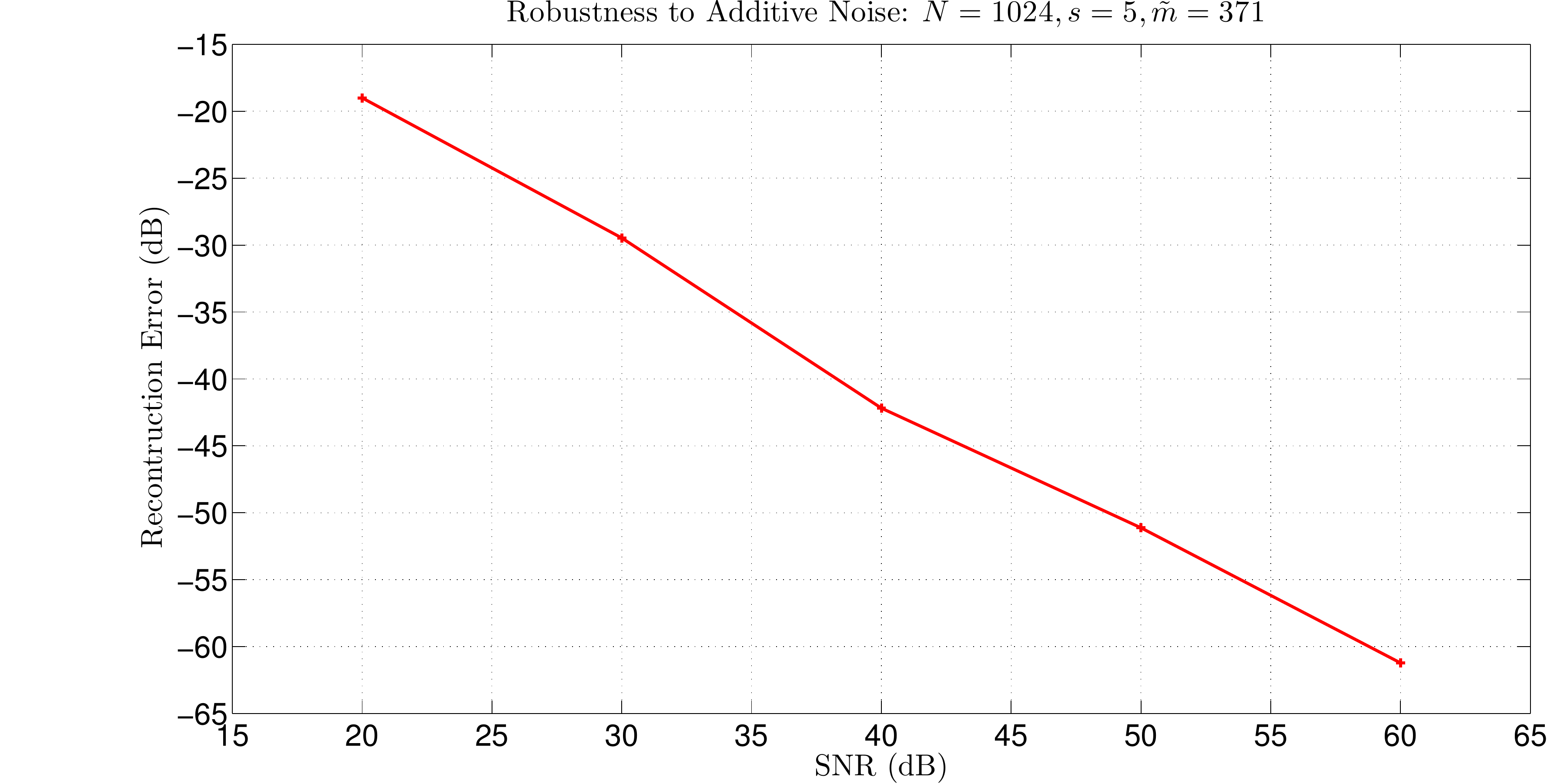}
\caption{Robustness to additive noise: $N = 1024, s = 5, \tilde m = \lceil 
	14s\log(N/s) \rceil$. }
\label{fig:robust}
\end{figure}

Figure \ref{fig:robust} illustrates the robustness of the recovery procedure to additive 
noise. We add i.i.d. zero-mean Gaussian noise at several signal-to-noise 
ratios (SNRs) to $\tilde m = \lceil 14s\log(N/s) \rceil$ magnitude measurements 
($N=1024, s=5, \tilde m=371$) and record relative reconstruction errors in decibels.
%
%$$	\mbox{Error (dB)} ~ = ~ 10 \log_{10} \left( || \tilde \x - \x ||_2^2 \right), $$
%where $\tilde \x$ denotes the recovered solution and $\x$ is the true signal. 
Each data point on the graph was obtained by averaging the results of 100 trials. We 
observe that the reconstruction error in every case is approximately equal to the 
added noise level, confirming the robust recovery properties of the proposed method.
\begin{figure}[hbtp]
        \centering
        \begin{subfigure}[b]{0.495\textwidth}
                \includegraphics[clip=true, trim = 0.5in 2.15in 0.25in 2.1in, scale=0.38]{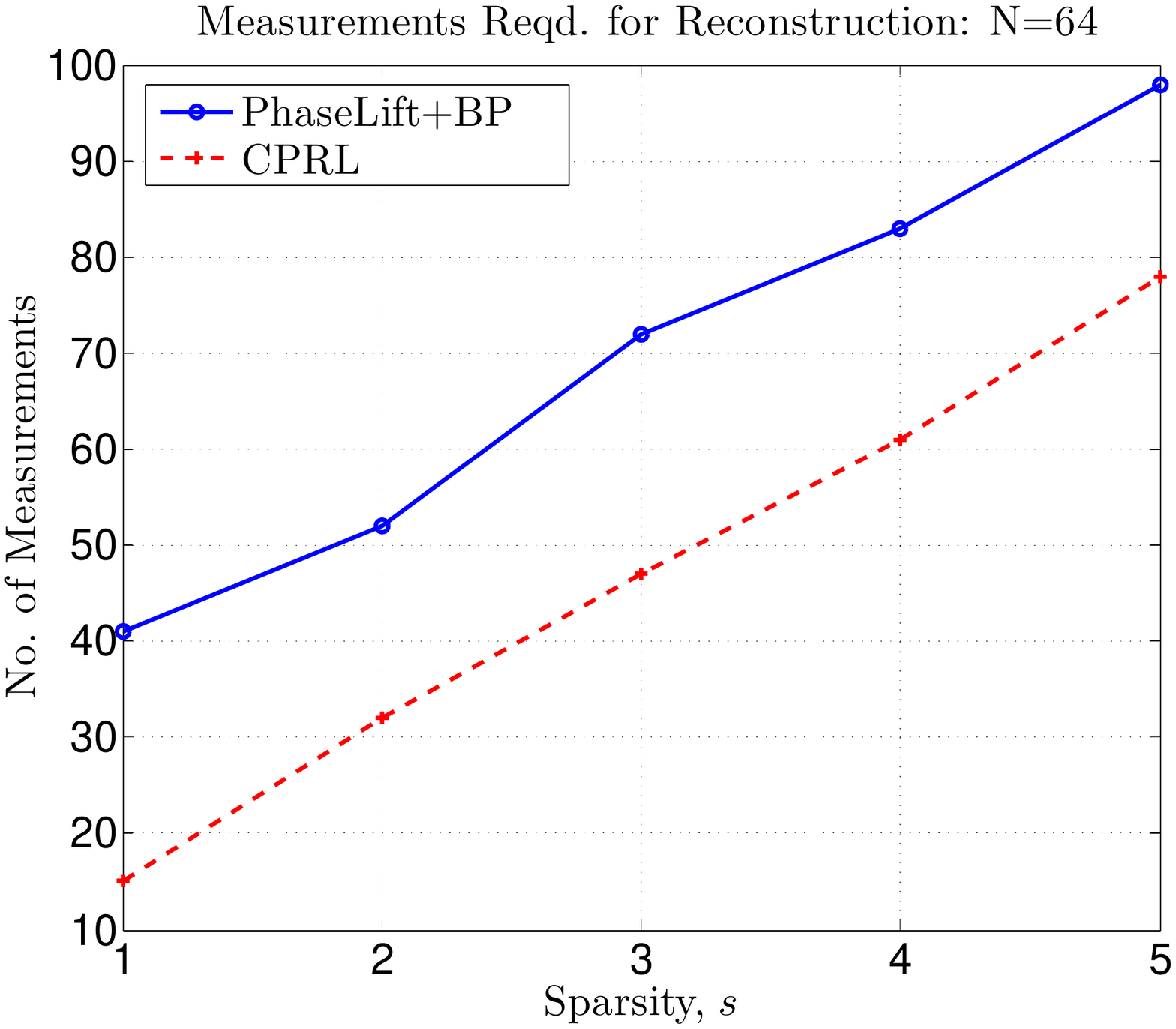}
                \caption{Minimum number of measurements}
                \label{fig:measurements}
        \end{subfigure}
        \hfill
        \begin{subfigure}[b]{0.495\textwidth}
                \includegraphics[clip=true, trim = 0in 2.15in 0.25in 2.1in, scale=0.38]{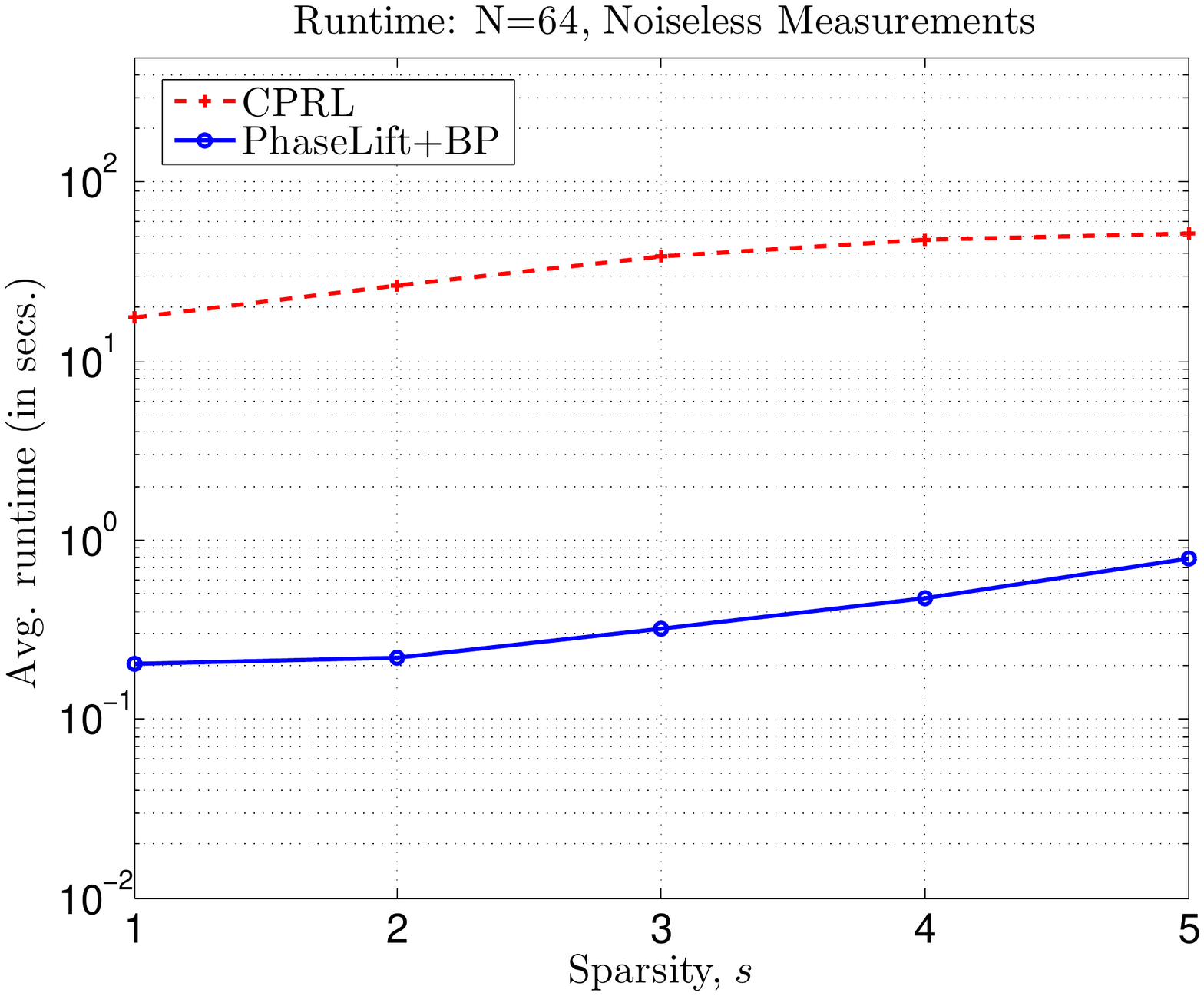}
                \caption{Runtime comparison}
                \label{fig:runtime}
        \end{subfigure}
        \caption{Runtime performance and minimum number of measurements required: 
	$N=64$, noiseless measurements}
        \label{fig:efficiency}
\end{figure}

Next, we demonstrate efficiency by plotting the average runtime and 
minimum number of measurements necessary for successful reconstruction. For the purposes of this 
discussion, we classify a reconstruction as successful if the relative $\ell_2$-norm error in the recovered signal is 
less than $10^{-5}$. We also provide comparisons with Compressive Phase Retrieval via Lifting ({\em CPRL}) 
\cite{ohlsson2012cprl}, an existing framework for sparse phase retrieval. Simulations were performed 
on a laptop computer with an Intel\textsuperscript\textregistered Core\textsuperscript{TM} i3-3120M processor, 
4GB RAM and Matlab R2014a. We first consider the reconstruction of an 
$s$-sparse signal ($N=64$) from perfect (noiseless) measurements. The minimum number of 
measurements\footnote{For the {\em PhaseLift}+BP implementation, we fixed the compressive 
sensing problem dimension to be $m~=~\lceil 1.75 s \log(N/s) \rceil$.} required for successful 
reconstruction is plotted in Figure \ref{fig:measurements}, while the corresponding runtime, 
averaged over 100 trials, is plotted in Figure \ref{fig:runtime}. Figure \ref{fig:measurements} was generated 
by starting with a small number of measurements, $\tilde m$, and incrementing this number to ensure successful 
reconstruction in at least 95 of the 100 trials. We notice that the {\em PhaseLift}+BP 
formulation requires a small number of additional measurements when compared to {\em CPRL}. This is potentially only 
the case for small values of $s$ since Theorem \ref{thm:CPR} shows that $\mathcal O(s \log (N/s))$ 
measurements suffice for the {\em PhaseLift}+BP formulation.  %, while {\em CPRL} requires $\mathcal O(s^2)$ 
%measurements. 
Moreover, since the {\em PhaseLift}+BP solution is obtained by solving a smaller SDP, the 
average runtime is significantly smaller (by several orders of magnitude) than {\em CPRL}, as 
shown in Figure \ref{fig:runtime}.

\section{Discussion}
\label{sec:Disscuss}

It is interesting to note that the compressive phase retrieval strategy discussed herein also immediately implies the existence of stable sublinear-time compressive phase retrieval algorithms.  These can be achieved by combining the phase retrieval technique of one's choice with a $o(N)$-time compressive sensing method (see, e.g., \cite{Iwen2014compressed}) in order to create a $o(N)$-time compressive phase retrieval algorithm.  In addition, we conclude by noting that random combinations of a random set of rows from a Fourier matrix will also exhibit the stable compressive phase retrieval property by Proposition~\ref{Prop}/Theorem~\ref{thm:CPR}.  This is of particular interest due to the special role that Fourier measurements play in many applications.

\bibliographystyle{abbrv}
\bibliography{SparsePR}

\end{document}